\documentclass[thmsa,11pt]{article}
\usepackage{amsmath}
\usepackage{amssymb}
\usepackage[letterpaper, margin=1.7in]{geometry}

\begin{document}

\title{ZFC proves that the class of ordinals is not weakly compact for
definable classes}
\author{Ali Enayat \& Joel David Hamkins}
\maketitle

\begin{abstract}
In $\mathrm{ZFC}$, the class $\mathrm{Ord}$ of ordinals is easily seen to
satisfy the definable version of strong inaccessibility. Here we explore
deeper $\mathrm{ZFC}$-verifiable combinatorial properties of $\mathrm{Ord}$,
as indicated in Theorems A \& B below. Note that Theorem A shows the
unexpected result that $\mathrm{Ord}$ is never definably weakly compact in
any model of $\mathrm{ZFC}$. \medskip

\noindent \textbf{Theorem A.~}\textit{Let }$\mathcal{M}$\textit{\ be any
model of }\textrm{ZFC}.

\medskip

\noindent \textbf{(1)} \textit{The definable tree property fails in }$%
\mathcal{M}$: \textit{There is an }$\mathcal{M}$\textit{-definable }$\mathrm{%
Ord}$-\textit{tree with no }$\mathcal{M}$-\textit{definable cofinal branch}.%
\textit{\medskip }

\noindent \textbf{(2) }\textit{The definable partition property fails in }$%
\mathcal{M}$: \textit{There is an }$\mathcal{M}$-\textit{definable }$2$%
\textit{-coloring }$f:[X]^{2}\rightarrow 2$\textit{\ for some }$\mathcal{M}$%
\textit{-definable proper class }$X$\textit{\ such that no }$\mathcal{M}$%
\textit{-definable proper classs is monochromatic for} $f$.\textit{\medskip }

\noindent \textbf{(3) }\textit{The definable compactness property for} $%
\mathcal{L}_{\infty \mathrm{,\omega }}$ \textit{fails in} $\mathcal{M}$:%
\textit{\ There is a definable theory }$\Gamma $\textit{\ in the logic} $%
\mathcal{L}_{\infty \mathrm{,\omega }}$ (\textit{in the sense of} $\mathcal{M%
}$) \textit{of size} $\mathrm{Ord}$ \textit{such that every set-sized
subtheory of} $\Gamma $ \textit{is satisfiable in }$\mathcal{M}$\textit{,
but there is no }$\mathcal{M}$\textit{-definable model of} $\Gamma $%
.\smallskip

\noindent \textbf{Theorem B.~}\textit{The definable} $\Diamond _{\mathrm{Ord}%
}$ \textit{principle holds in a model }$\mathcal{M}$\textit{\ of }\textrm{ZFC%
}\textit{\ iff }$\mathcal{M}$\textit{\ carries an }$\mathcal{M}$\textit{%
-definable} \textit{global well-ordering}.

\medskip

Theorems A and B above can be recast as theorem schemes in $\mathrm{ZFC}$,
or as asserting that a single statement in the language of class theory
holds in all `spartan' models of $\mathrm{GB}$ (G\"{o}del-Bernays class
theory); where a spartan model of $\mathrm{GB}$\ is any structure of the
form $(\mathcal{M},D_{\mathcal{M}})$\textit{,} where $\mathcal{M}\models
\mathrm{ZF}$\textit{\ }and $D_{\mathcal{M}}$ is the family of\textit{\ }$%
\mathcal{M}$-definable classes. Theorem C gauges the complexity of the
collection $\mathrm{GB}_{\mathrm{spa}}$ of (G\"{o}del-numbers of) sentences
that hold in all spartan models of $\mathrm{GB.}$\textit{\medskip }

\noindent \textbf{Theorem C.~}$\mathrm{GB}_{\mathrm{spa}}$ \textit{is} $\Pi
_{1}^{1}$\textit{-complete.}

\let\thefootnote\relax\footnotetext{\textsc{%
}}

\let\thefootnote\relax

\end{abstract}

\pagebreak

\textbf{1.~Introduction \& Preliminaries}\bigskip

In $\mathrm{ZFC}$, the class $\mathrm{Ord}$ of ordinals satisfies the
definable version of strong inaccessibility since the power set axiom and
the axiom of choice together make it evident that $\mathrm{Ord}$ is closed
under cardinal exponentiation; and the scheme of replacement ensures the
definable regularity of $\mathrm{Ord}$\textbf{\ }in the sense that for each
cardinal $\kappa <\mathrm{Ord}$, the range of every definable ordinal-valued
map $f$ with domain $\kappa $ is bounded in $\mathrm{Ord}$. In this paper we
investigate more subtle definable combinatorial properties of $\mathrm{Ord}$%
\textbf{\ }in the context of $\mathrm{ZFC}$ to obtain results, each of which
takes the form of a \textit{theorem scheme} within $\mathrm{ZFC}$. In
Section 2 we establish a number of results that culminate in Theorem 2.6,
which states that the tree property fails for definable classes across all
models of $\mathrm{ZFC}$; this result is then used in Section 3 to show the
failure of the partition property for definable classes, and the failure of
weak compactness of $\mathrm{Ord}$ for definable classes in all models of $%
\mathrm{ZFC.}$ \emph{Thus, the results in Sections 2 and 3 together
demonstrate the unexpected }$\mathrm{ZFC}$\emph{-provable failure of the
definable version of a large cardinal property for} $\mathrm{Ord}$. In
Section 4 we establish the equivalence of the combinatorial principle $%
\Diamond _{\mathrm{Ord}}$ and the existence of a definable global choice
function across all models of $\mathrm{ZFC.}$ \medskip

The results in Sections 2 through 4 can be viewed as stating that certain
sentences in the language of class theory hold in all `spartan' models of
\textrm{GB} (G\"{o}del-Bernays class theory), i.e., in all models of $%
\mathrm{GB}$ of the form $(\mathcal{M},\mathcal{D}_{\mathcal{M}})$, where $%
\mathcal{M}$ is a model of $\mathrm{ZF}$ and $\mathcal{D}_{\mathcal{M}}$ is
the collection of $\mathcal{M}$-definable subsets of $M$. For example
Theorem 2.6 is equivalent to the veracity of the statement \textquotedblleft
if the axiom of choice for sets holds, then there is an $\mathrm{Ord}$%
-Aronszajn tree\textquotedblright\ in every spartan model of \textrm{GB. }In
Section 5 we show that the \textit{theory} of all spartan models of $\mathrm{%
GB}$, when viewed as a subset of $\omega $ via G\"{o}del-numbering, is $\Pi
_{1}^{1}$-complete; and a fortiori, it is not computably axiomatizable.
\medskip

We now turn to reviewing pertinent preliminaries concerning models of set
theory. Our meta-theory is $\mathrm{ZFC}$.\medskip

\noindent \textbf{1.1.~Definition.~}Suppose $\mathcal{M}=(M,\in ^{\mathcal{M}%
})$\ and $\mathcal{N}=(N,\in ^{\mathcal{N}})$ are models of set theory. Note
that we are not assuming that either $\mathcal{M}$ or $\mathcal{N}$\ is
well-founded.\textbf{\medskip }

\noindent \textbf{(a)} For $m\in M,$ let $m_{\mathcal{M}}:=\{x\in M:x\in ^{%
\mathcal{M}}m\}.$ If $\mathcal{M}$ $\subseteq \mathcal{N}$ (i.e., $\mathcal{M%
}$ is a submodel of $\mathcal{N}$) and $m\in M$, then $\mathcal{N}$ \textit{%
fixes} $m$ if $m_{\mathcal{M}}=m_{\mathcal{N}}.$ $\mathcal{N}$ \textit{end
extends} $\mathcal{M}$, written $\mathcal{M}\subseteq _{e}\mathcal{N}$, iff $%
\mathcal{N}$ fixes every $m\in M.$ Equivalently: $\mathcal{M}\subseteq _{e}%
\mathcal{N}$ iff $\mathcal{M}$ is a transitive submodel of $\mathcal{N}$ in
the sense that if $x\in ^{\mathcal{N}}y$ for some $x\in N$ and some $y\in M,$
then $x\in ^{\mathcal{M}}y$. \textbf{\medskip }

\noindent \textbf{(b)} Given $n\in \omega $, $\mathcal{N}$ is a \textit{%
proper }$\Sigma _{n}$\textit{-e.e.e.}\textbf{~}of $\mathcal{M}$
(\textquotedblleft e.e.e.\textquotedblright\ stands for \textquotedblleft
elementary end extension\textquotedblright ), iff $\mathcal{M\subsetneq }_{e}%
\mathcal{N}$, and $\mathcal{M}\prec _{\Sigma _{n}}\mathcal{N}$ (i.e., $%
\Sigma _{n}$-statements with parameters from $M$ are absolute in the passage
between $\mathcal{N}$ and $\mathcal{M)}$. It is well-known that if $\mathcal{%
M}\prec _{\Sigma _{2}}\mathcal{N}$ and $\mathcal{M}\models \mathrm{ZF}$,
then $\mathcal{N}$ is a \textit{rank extension} of $\mathcal{M}$, i.e.,
whenever $a\in M$ and $b\in N\backslash M$, then $\mathcal{N}\models \rho
(a)\in \rho (b),$ where $\rho $ is the usual ordinal-valued rank function on
sets.\textbf{\medskip }

\noindent \textbf{(c)} Given $\alpha \in \mathrm{Ord}^{\mathcal{M}}$, $%
\mathcal{M}_{\alpha }$ denotes the structure $(\mathrm{V}_{\alpha },\in )^{%
\mathcal{M}},$ and $M_{\alpha }=\mathrm{V}_{\alpha }^{M}$.\textbf{\medskip }

\noindent \textbf{(d)} For $X\subseteq M^{n}$ (where $n\in \omega ),$ we say
that $X$ is $\mathcal{M}$-\textit{definable} iff $X$\ is parametrically
definable in $\mathcal{M}$.\textbf{\medskip }

\noindent \textbf{(e)} $\mathcal{N}$ is a\textit{\ conservative }extension
of $\mathcal{M}$, written $\mathcal{M}\subseteq _{\mathrm{cons}}\mathcal{N}$%
, iff the intersection of any $\mathcal{N}$-definable subset of $N$ with $M$%
\emph{\ }is $\mathcal{M}$-definable.\medskip

\noindent For models of $\mathrm{ZF}$, the set-theoretical sentence $\exists
p\left( \mathrm{V}=\mathrm{HOD}(p)\right) $ expresses: \textquotedblleft
there is some $p$ such that every set is first order definable in some
structure of the form $(\mathrm{V}_{\alpha },\in ,p)$ with $p\in \mathrm{V}%
_{\alpha }$\textquotedblright . The following theorem is well-known; the
equivalence of (a) and (b) will be revisited in Theorem 4.2.\medskip

\noindent \textbf{1.2.~Theorem.~}\textit{The following statements are
equivalent for} $\mathcal{M}\models \mathrm{ZF}$\textrm{:}\textbf{\medskip }

\noindent \textbf{(a) }$\mathcal{M}\models \exists p\left( \mathrm{V}=%
\mathrm{HOD}(p)\right) .$\textbf{\medskip }

\noindent \textbf{(b) }\textit{For some} $p\in M$ \textit{and some
set-theoretic formula} $\varphi (x,y,\overline{p})$ (\textit{where} $%
\overline{p}$ \textit{is a name for} $p$) $\mathcal{M}$ \textit{satisfies }%
\textquotedblleft $\varphi $ \textit{well-orders the universe}%
\textquotedblright .\textbf{\medskip }

\noindent \textbf{(c)} \textit{For some} $p\in M$ \textit{and some} $\Sigma
_{2}$-\textit{formula} $\varphi (x,y,\overline{p})$ $\mathcal{M}$ \textit{%
satisfies }\textquotedblleft $\varphi $ \textit{well-orders the universe}%
\textquotedblright .\textbf{\medskip }

\noindent (\textbf{d)} $\mathcal{M}\models \forall x(x\neq \varnothing
\rightarrow f(x)\in x)$ \textit{for some }$\mathcal{M}$-\textit{definable} $%
f:M\rightarrow M.$\medskip

\noindent Next we use definable classes to lift certain combinatorial
properties of cardinals to the class of ordinals. \medskip

\noindent \textbf{1.3.~Definitions.~}Suppose\textit{\ }$\mathcal{M}\models
\mathrm{ZFC.}$\medskip

\noindent \textbf{(a)} Suppose $\tau =\left( T,\ <_{T}\right) $ is a tree
ordering, where both $T$ and $<_{T}$ are $\mathcal{M}$-definable. $\tau $ is
an $\mathrm{Ord}$-\textit{tree in }$\mathcal{M}$ iff $\mathcal{M}$ satisfies
\textquotedblleft $\tau $ is a well-founded tree of height $\mathrm{Ord}$%
\textbf{\ }and for all $\alpha \in \mathrm{Ord}$, the collection $T_{\alpha
} $ of elements of $T$ at level $\alpha $ of $\tau $ form a
set\textquotedblright . Such a tree $\tau $ is said to be a \textit{definably%
} $\mathrm{Ord}$-\textit{Aronszajn} tree in $\mathcal{M}$ iff no cofinal
branch of $\tau $ is $\mathcal{M}$-definable.\medskip

\noindent \textbf{(b)} \textit{The definable tree property for }$\mathrm{Ord}
$ \textit{fails} \textit{in} $\mathcal{M}$\textit{\ }iff there exists a
definably $\mathrm{Ord}$-\textit{Aronszajn} tree in $\mathcal{M}$.\footnote{%
This notion should not be confused with the \textit{definable tree property
of a cardinal} $\kappa $, first introduced and studied by Leshem $\cite%
{Leshem},$ which stipulates that every $\kappa $-tree that is first order
definable (parameters allowed) in the structure $\left( H(\kappa ),\in
\right) $ has a cofinal branch $B$ (where $H(\kappa )$ is the collection of
sets that are hereditarily of cardinality less than $\kappa $). Note that in
this definition $B$ is not required to be first order definable in $\left(
H(\kappa ),\in \right) ;$ so every weakly compact cardinal has the definable
tree property.
\par
{}}\medskip

\noindent \textbf{(c) }\textit{The definable proper class partition property
fails in }$\mathcal{M}$ iff there is an $\mathcal{M}$-definable proper class
$X$ of $M$ with an $\mathcal{M}$-definable $2$-coloring $f:[X]^{2}%
\rightarrow 2$ such that there is no $\mathcal{M}$-definable monochromatic
proper class for $f$. We also say that $\mathrm{Ord}\rightarrow \left(
\mathrm{Ord}\right) _{2}^{2}$ \textit{fails in} $\mathcal{M}$ iff there is
an $\mathcal{M}$-definable $2$-coloring $f:[\mathrm{Ord}]^{2}\rightarrow 2$
such that there is no $\mathcal{M}$-definable monochromatic proper class for
$f$. \medskip

\noindent \textbf{(d) }\textit{The definable} \textit{compactness property
for }$\mathcal{L}_{\mathrm{\infty },\mathrm{\omega }}$ \textit{fails in }$%
\mathcal{M}$ iff there is an\textit{\ }$\mathcal{M}$-definable theory $%
\Gamma $ formulated in the logic $\mathcal{L}_{\mathrm{\infty },\mathrm{%
\omega }}$ (in the sense of $\mathcal{M}$) such that every set-sized
subtheory of $\Gamma $ is satisfiable in $\mathcal{M}$, but there is no%
\textit{\ }$\mathcal{M}$\textit{-}definable model\textit{\ }of $T$. Here $%
\mathcal{L}_{\mathrm{\infty },\mathrm{\omega }}$ is the extension of first
order logic that allows conjunctions and disjunctions applied to \textit{sets%
} of formulae (of any cardinality) with only a finite number of free
variables, as in \cite[Ch.III]{Barwise}. \medskip

\noindent \textbf{(e) }An $\mathcal{M}$-definable subset $E$ of $\mathrm{Ord}%
^{\mathcal{M}}$ is said to be \textit{definably} $\mathcal{M}$-\textit{%
stationary} iff $E\cap C\neq \varnothing $ for every $\mathcal{M}$-definable
subset $C$ of $\mathcal{M}$ such that $C$ is closed and unbounded in $%
\mathrm{Ord}^{\mathcal{M}}$.\medskip

\noindent \textbf{(f) }\textit{The definable}\textbf{\ }$\Diamond _{\mathrm{%
Ord}}$ \textit{holds in }$\mathcal{M}$ iff there is some $\mathcal{M}$%
-definable $\vec{A}=\left\langle A_{\alpha }:\alpha \in \mathrm{Ord}^{%
\mathcal{M}}\right\rangle $ such that $\mathcal{M}$ satisfies
\textquotedblleft $A_{\alpha }\subseteq \alpha $ for all $\alpha \in \mathrm{%
Ord}$\textquotedblright , and for all $\mathcal{M}$-definable $A\subseteq
\mathrm{Ord}^{\mathcal{M}}$ there is $E\subseteq \mathrm{Ord}^{\mathcal{M}}$
such that $E$ is definably $\mathcal{M}$-stationary and $A_{\alpha }=A\cap
\alpha $ for all $\alpha \in E$. Here $\vec{A}$ is said to be $\mathcal{M}$%
-definable if there is an $\mathcal{M}$-definable $A$ such that $A_{\alpha
}=\{m:\left\langle m,\alpha \right\rangle \in A\}$ for each $\alpha \in
\mathrm{Ord}^{\mathcal{M}}.$ \medskip

\bigskip

\begin{center}
\textbf{2.~The failure of the definable tree property for the class of
ordinals}\bigskip
\end{center}

The proof of the main result of this section (Theorem 2.6) is based on a
number of preliminary model-theoretic results which are of interest in their
own right. We should point out that a proof of a special case of Theorem 2.6
was sketched in \cite[Remark 3.5]{Ali Power-like} for models of set theory
with built-in global choice functions, using a more technical argument than
the one presented here.\medskip

We begin with the following theorem which refines a result of Kaufmann \cite[%
Theorem 4.6]{Matt Topless and blunt}. The proof uses an adaptation of
Kaufmann's proof based on a strategy introduced in \cite[Theorem 1.5(a)]{Ali
TAMS}\textbf{. }\medskip

\noindent \textbf{2.1.~Theorem.~}\textit{No model of }$\mathrm{ZFC}$ \textit{%
has a proper conservative} $\Sigma _{3}$-\textit{e.e.e}.\textit{\ }\medskip

\noindent \textbf{Proof. }Suppose to the contrary that $\mathcal{M}\models
\mathrm{ZF}$ and $\mathcal{M}\prec _{\Sigma _{3},\mathrm{e},\mathrm{cons}}%
\mathcal{N}$ for some $\mathcal{N}$. Let $\varphi $ be the statement that
expresses the following instance of the reflection theorem:

\begin{center}
$\forall \lambda \in \mathrm{Ord}\mathbf{\ }\exists \beta \in \mathrm{Ord}%
\left( \lambda \in \beta \wedge \ (\mathrm{V}_{\beta },\in )\prec _{\Sigma
_{1}}\left( \mathrm{V},\in \right) \right) .$
\end{center}

\noindent Using the fact that the satisfaction predicate for $\Sigma _{1}$%
-formulae is $\Sigma _{1}$-definable it is easy to see that $\varphi $ is a $%
\Pi _{3}$-statement, and thus $\varphi $ also holds in $\mathcal{N}$ since $%
\varphi $ holds in $\mathcal{M}$ by the reflection theorem.\footnote{%
Recall that, provably in \textrm{ZF}, the ordinals $\beta $ such that $(%
\mathrm{V}_{\beta },\in )\prec _{\Sigma _{1}}\left( \mathrm{V},\in \right) $
are precisely the fixed points of the $\beth $-function.} So we can fix some
$\lambda \in \mathrm{Ord}^{\mathcal{N}}\backslash \mathrm{Ord}^{\mathcal{M}}$
and some $\mathcal{N}$-ordinal $\beta >\lambda $ of $\mathcal{N}$ such that:

\begin{center}
$\mathcal{N}_{\beta }\prec _{\Sigma _{1}}\mathcal{N}$.
\end{center}

\noindent Note that this implies that $\mathcal{N}_{\beta }$ can
meaningfully define the satisfaction predicate for every set-structure
`living in'\ $\mathcal{N}_{\beta }$ since that $\mathcal{N}_{\beta }$ is a
model of a substantial fragment of $\mathrm{ZF}$, including $\mathrm{KP}$
(Kripke-Platek set theory), and already $\mathrm{KP}$\ is sufficient for
this purpose \cite[III.2]{Barwise}. Also, since the statement
\textquotedblleft every set can be well-ordered\textquotedblright\ is a $\Pi
_{2}$-statement which holds in $\mathcal{M}$ by assumption, it also holds in
$\mathcal{N}$, and therefore we can fix a binary relation $w$ in $\mathcal{N}
$ such that, as viewed in $\mathcal{N}$, $w$ is a well-ordering of $\mathrm{V%
}_{\beta }.$ Hence for any $\alpha \in \mathrm{Ord}^{\mathcal{M}}$ with $%
\alpha <\beta $, \textit{within} $\mathcal{N}$ one can define the submodel $%
\mathcal{K}_{\alpha }$ of $\mathcal{N}_{\beta }$ whose universe $K_{\alpha }$
is defined via:

\begin{center}
$K_{\alpha }:=\{a\in \mathrm{V}_{\beta }:a$ is first order definable in $(%
\mathcal{N}_{\beta },w,\lambda ,m)_{m\in \mathrm{V}_{\alpha }}$\}.
\end{center}

\noindent Clearly $M_{\alpha }\cup \{\lambda \}\subsetneq \mathcal{K}%
_{\alpha }\prec \mathcal{N}_{\beta }$, and of course $\mathcal{K}_{\alpha }$
is a member of $\mathcal{N}$. Next let:

\begin{center}
$\mathcal{K}:=\bigcup\limits_{\alpha \in \mathrm{Ord}^{\mathcal{M}}}\mathcal{%
K}_{\alpha }.$
\end{center}

\noindent Note that we have:

\begin{center}
$\mathcal{M}\subsetneq _{\text{e}}\mathcal{K}\preceq \mathcal{N}_{\beta
}\prec _{1}\mathcal{N}.$
\end{center}

\noindent We now make a crucial case distinction: either (a) $\mathrm{Ord}^{%
\mathcal{K}}\backslash \mathrm{Ord}^{\mathcal{M}}$ has minimum element, or
(b) it does not. The proof will be complete once we verify that both cases
lead to a contradiction.\medskip

\textbf{Case (a)}. Let $\eta =\min (\mathrm{Ord}^{\mathcal{K}}\backslash
\mathrm{Ord}^{\mathcal{M}})$. We claim that $\mathcal{M\prec N}_{\eta }$. To
see this, we use Tarski's test for elementarity: suppose $\mathcal{N}_{\eta
}\models \exists x\varphi (x,\overline{m})$ for some $m\in M$ and some
formula $\varphi (x,y),$ and let $\theta _{0}$ be defined in $\mathcal{N}%
_{\beta }$ as the least ordinal $\theta $ such that $x\in \mathrm{V}_{\theta
}$ and $\mathrm{V}_{\eta }\models \exists x$ $\varphi (x,\overline{m})$.
Then $\theta _{0}\in K$ and clearly $\theta _{0}<\eta $, which shows that $%
\theta _{0}\in \mathrm{Ord}^{\mathcal{M}}$. Hence $\mathcal{N}_{\eta
}\models \varphi (\overline{m_{0}},\overline{m})$ for some $m_{0}\in M,$
thus completing the proof of $\mathcal{M\prec N}_{\eta }.$ But if $\mathcal{%
M\prec N}_{\eta }$, then we can choose $S$ in $\mathcal{N}$ such that:

\begin{center}
$\mathcal{N}\models S=\{\ulcorner \varphi (\overline{m})\urcorner \in
\mathrm{V}_{\eta }:\mathcal{N}\models $ \textquotedblleft $(\mathrm{V}_{\eta
},\in )\models \varphi (\overline{m})$\textquotedblright $\}.$
\end{center}

\noindent Based on the assumption that $\mathcal{N}$\ is a conservative
extension of $\mathcal{M}$, $S\cap M$ should be an $\mathcal{M}$-definable
satisfaction predicate for $\mathcal{M}$, which contradicts (a version of)
Tarski's undefinability of truth theorem. \medskip

\textbf{Case (b)}. This is the more difficult case, where $\mathrm{Ord}^{%
\mathcal{K}}\backslash \mathrm{Ord}^{\mathcal{M}}$ has no least element. Let
$\Phi :=\bigcup\limits_{\alpha \in \mathrm{Ord}^{\mathcal{M}}}\Phi _{\alpha
} $, where

\begin{center}
$\Phi _{\alpha }:=\{\ulcorner \varphi (c,\overline{m})\urcorner \in M:%
\mathcal{N}\models $ \textquotedblleft $(\mathrm{V}_{\beta },\in ,w,\lambda
,m)_{m\in \mathrm{V}_{\alpha }}\models \varphi (c,\overline{m})$%
\textquotedblright $\}.$
\end{center}

\noindent In the above definition of $\Phi _{\alpha }$, the constant $c$ is
interpreted as $\lambda $ and $\varphi (c,\overline{m})$ ranges over first
order formulae in the sense of $\mathcal{M}$ (or equivalently: in the sense
of $\mathcal{N}$) in the language

\begin{center}
$\mathcal{L}_{\alpha }=\{\in ,\vartriangleleft ,c\}\cup \{\overline{m}:m\in
\mathrm{V}_{\alpha }\}$,
\end{center}

\noindent where $c$ is a new constant symbol and $\vartriangleleft $ is a
binary relation symbol interpreted by $w$. Thus $\Phi $ can be thought of as
the type of $\lambda $ in $\mathcal{N}_{\beta }$ over $M$. Since $\mathcal{N}
$ is assumed to be a conservative extension of $\mathcal{M}$, $\Phi $ is $%
\mathcal{M}$-definable via some unary formula $\phi $. Hence $\Gamma $ below
is also $\mathcal{M}$-definable via some unary formula $\gamma :$

\begin{center}
$\overset{\Gamma }{\overbrace{\left\{ \ulcorner t(c,\overline{m})\urcorner
\in M:\phi \left( \ulcorner t(c,\overline{m})\in \mathrm{Ord}\mathbf{%
\urcorner }\right) \ \mathrm{and}\ \forall \theta \in \mathrm{Ord}(\phi
\left( \ulcorner t(c,\overline{m})>\overline{\theta }\mathbf{\urcorner }%
\right) \right\} }},$
\end{center}

\noindent where $t$ is a definable term in the language $\mathcal{L}$, i.e.,
$t(c,\overline{m})$ is an $\mathcal{L}$-definition $\varphi (c,\overline{m}%
,x)$ of some element $x$. So, officially speaking, $\Gamma $ consists of $%
\ulcorner \varphi (c,\overline{m},x)\urcorner \in M$ that satisfy the
following three conditions:\medskip

\noindent (1) $\phi \left( \ulcorner \exists !x\varphi (c,\overline{m},x%
\mathbf{\urcorner }\right) .$

\noindent (2) $\phi \left( \ulcorner \forall x\left( \varphi (c,\overline{m}%
,x)\rightarrow x\in \mathrm{Ord}\right) \mathbf{\urcorner }\right) .$

\noindent (3) $\forall \theta \in \mathrm{Ord}\ \phi \left( \ulcorner
\forall x\left( \varphi (c,\overline{m},x)\rightarrow x>\overline{\theta }%
\right) \mathbf{\urcorner }\right) .$\medskip

\noindent Since $\mathrm{Ord}^{\mathcal{K}}\backslash \mathrm{Ord}^{\mathcal{%
M}}$ has no minimum element (recall: we are analysing case (b)), $\mathcal{M}%
\models \psi $, where:

\begin{center}
$\psi :=\forall t\left( \gamma (t)\rightarrow \exists t^{\prime }(\gamma
(t^{\prime })\wedge \phi (\ulcorner t^{\prime }\in t\mathbf{\urcorner }%
)\right) .$
\end{center}

\noindent Choose $k$ such that $\psi $ is a $\Sigma _{k}$-statement, and use
the reflection theorem in $\mathcal{M}$ to pick $\mu \in \mathrm{Ord}^{%
\mathcal{M}}$ such that $\mathcal{M}_{\mu }\prec _{\Sigma _{k}}\mathcal{M}.$
Then $\psi $ holds in $\mathcal{M}_{\mu }$, so by DC (dependent choice,
which holds in $\mathcal{M}$ since AC holds in $\mathcal{M}$), there is some
function $f_{c}$ in $\mathcal{M}$ such that:

\begin{center}
$\mathcal{M}\models \forall n\in \omega \ \phi \left( \ulcorner
f_{c}(n+1)\in f_{c}(n)\mathbf{\urcorner }\right) .$
\end{center}

\noindent Let $\alpha \in \mathrm{Ord}^{\mathcal{M}}$ be large enough so
that $M_{\alpha }$ contains all constants $\overline{m}$ that occur in any
of the terms in the range of $f$; let $f_{\lambda }(n)$ be defined in $%
\mathcal{N}$ as the result of replacing all occurrences of the constant $c$
with $\overline{\lambda }$ in $f_{c}(n)$; and let $g(n)$ be defined in $%
\mathcal{N}$ as the interpretation of $f_{\lambda }(n)$ in $(\mathrm{V}%
_{\beta },\in ,w,\lambda ,m)_{m\in \mathrm{V}_{\alpha }}.$ Then $\mathcal{N}$
satisfies:

\begin{center}
$\forall n\in \omega \ \left( g(n)\in g(n+1)\right) $,
\end{center}

\noindent which contradicts the foundation axiom in $\mathcal{N}$. The proof
is now complete.\hfill $\square $\bigskip

\noindent \textbf{2.2.~Definition.}

\noindent \textbf{(a)} Given ordinals $\alpha <\beta ,$ $\mathcal{V}_{\beta
,\alpha }$ denotes the structure $(\mathrm{V}_{\beta },\in ,a)_{a\in \mathrm{%
V}_{\alpha }}$, and for a model $\mathcal{M}\models \mathrm{ZF},$

\begin{center}
$\mathcal{M}_{\beta ,\alpha }:=\left( \mathcal{V}_{\beta ,\alpha }\right) ^{%
\mathcal{M}}.$
\end{center}

\noindent \textbf{(b)} Given a meta-theoretic natural number $n$, $\tau _{n}$
denotes the definable tree whose nodes at level $\alpha $ consist of first
order theories of the form $\mathrm{Th}(\mathcal{V}_{\beta ,\alpha },s)$,
where $s\in \mathrm{V}_{\beta }\backslash \mathrm{V}_{\alpha },$ and $\beta $
is $n$-correct\footnote{%
An ordinal $\beta $ is $n$-correct when $\left( \mathrm{V}_{\beta },\in
\right) \prec _{\Sigma _{n}}\left( \mathrm{V},\in \right) .$}. The language
of $\mathrm{Th}(\mathcal{V}_{\beta ,\alpha },s)$ consists of $\{\in \}$ plus
constants $\overline{m}$ for each $m\in \mathrm{V}_{\alpha },$ and a new
constant $c$ whose denotation is $s$. The ordering of the tree is by
set-inclusion.\smallskip

\noindent \textbf{2.3.~Lemma.~}\textit{For each meta-theoretic natural
number }$n$, $\mathrm{ZFC}$ \textit{proves} \textquotedblleft $\tau _{n}$
\textit{is an} $\mathrm{Ord}$-\textit{tree}\textquotedblright .\medskip

\noindent \textbf{Proof.~}Thanks to the Montague-Vaught reflection theorem,
there are plenty of nodes at any ordinal level $\alpha $. On the other hand,
since each $\mathrm{Th}(\mathcal{V}_{\beta ,\alpha },s)$ can be canonically
coded as a subset of $\mathrm{V}_{\alpha },$ and $\left\vert \mathrm{V}%
_{\omega +\alpha }\right\vert =\beth _{\alpha },$ there are at most $\beth
_{\alpha }$-many nodes at level $\alpha $ \hfill $\square $\medskip

\noindent \textbf{2.4}.\textbf{~Remark.~}One may `prune' every $\mathrm{Ord}$%
-tree $\tau $ to obtain a definable subtree $\tau ^{\ast }$ which has nodes
of arbitrarily high level in $\mathrm{Ord}$ by simply throwing away the
nodes whose set of successors have bounded height and then using the
replacement scheme to verify that the subtree $\tau ^{\ast }$ thus obtained
has height $\mathrm{Ord}$. See \cite[Lemma 3.11]{Kunen Text} for a similar
construction for $\kappa $-trees (where $\kappa $ is a regular cardinal).
\medskip

\noindent \textbf{2.5.~Lemma.~}\textit{Suppose} $\mathcal{M}$ \textit{is a
model of} $\mathrm{ZFC}$ \textit{that carries an }$\mathcal{M}$\textit{%
-definable global well-ordering}. \textit{Furthermore, suppose that }$n\geq
3 $ \textit{and the tree }$\tau _{n}^{\mathcal{M}}$ \textit{has a branch }$B$%
.\textit{\ \textbf{Then:\smallskip }}

\noindent \textbf{(a)}\textit{\ There is a model }$\mathcal{N}$\textit{\ and
a proper embedding }$j:\mathcal{M}\rightarrow \mathcal{N}$ \textit{such that}
$j(\mathcal{M})\prec _{\mathrm{e},n}\mathcal{N}$\textit{.\textbf{\smallskip }%
}

\noindent \textbf{(b)} \textit{Both }$\mathcal{N}$\textit{\ and }$j$ \textit{%
are} $\mathcal{M}$-\textit{definable if }$B$\textit{\ is }$\mathcal{M}$%
\textit{-definable.\textbf{\smallskip }}

\noindent \textbf{(c)} $\mathcal{N}$\textit{\ is a conservative extension of
}$j(\mathcal{M})$\textit{\ if }$B$\textit{\ is }$\mathcal{M}$\textit{%
-definable.}\medskip

\noindent \textbf{Proof.~}We will only prove (a) since the proof of (b) will
be clear by an inspection of the proof of (a), and (c) is an immediate
consequence of (b). Let $B$ be a branch of $\tau _{n}^{\mathcal{M}}.$ Each
node in $B$ is a first order theory in the sense of $\mathcal{M}$ and is of
the form $\left( \mathrm{Th}(\mathcal{\mathrm{V}}_{\beta ,\alpha },s)\right)
^{\mathcal{M}}$. Note that $\left( \mathrm{Th}(\mathcal{V}_{\beta ,\alpha
},s)\right) ^{\mathcal{M}}$ is not the necessarily the same as $\mathrm{Th}(%
\mathcal{M}_{\beta ,\alpha },s)$, since the latter is the collection of
\textit{standard} sentences in $\left( \mathrm{Th}(\mathcal{V}_{\beta
,\alpha },s)\right) ^{\mathcal{M}}$. In particular, if $\mathcal{M}$ not $%
\omega $-standard, then:

\begin{center}
$\mathrm{Th}(\mathcal{M}_{\beta ,\alpha },s)\subsetneq \left( \mathrm{Th}(%
\mathcal{V}_{\beta ,\alpha },s)\right) ^{\mathcal{M}}.$
\end{center}

\noindent For each $\alpha \in \mathrm{Ord}^{\mathcal{M}}$, let $b_{\alpha }$
be the node of $B$ at level $\alpha $. We may choose some $\beta _{\alpha
}\in \mathrm{Ord}^{\mathcal{M}}$ and some $s_{\alpha }\in \left( \mathrm{V}%
_{\beta _{\alpha }}\backslash \mathrm{V}_{\alpha }\right) ^{\mathcal{M}}$
such that:

\begin{center}
$b_{\alpha }=\left( \mathrm{Th}(\mathcal{V}_{\beta _{\alpha },\alpha
},s_{\alpha })\right) ^{\mathcal{M}}.$
\end{center}

\noindent The above choices of $\beta _{\alpha }$ and $s_{\alpha }$ are
performed at the meta-theoretic level (where $\mathrm{ZFC}$\ is assumed);
however if $B$ is $\mathcal{M}$-definable, then so is the map $\alpha
\mapsto b_{\alpha }$, which in turn shows that the maps $\alpha \mapsto
\beta _{\alpha }$ and $\alpha \mapsto s_{\alpha }$ can also be arranged to
be $\mathcal{M}$-definable since $\mathcal{M}$ is assumed to carry an $%
\mathcal{M}$-definable global well-ordering (the definability of these two
maps plays a key role in verifying that an inspection of the proof of (a)
yields a proof of (b)).\smallskip

We now explain how to use $B$ to construct the desired structure $\mathcal{N}
$. In order to do so, we need some definitions:\medskip

\noindent $(i)$ Let $\mathcal{L}$ be the language consisting of the usual
language $\{\in \}$ of set theory, augmented with a binary relation symbol $%
\vartriangleleft $, constants $\overline{m}$ for each $m\in M,$ and a new
constant $c.$ \medskip

\noindent $(ii)$ For each $\alpha \in \mathrm{Ord}^{\mathcal{M}}$ let $%
\mathcal{N}_{\alpha }$ be the submodel of $\mathcal{M}_{\beta _{\alpha }}$
whose universe $N_{\alpha }$ consists of elements of $M_{\beta _{\alpha }}$
that are first order definable in the structure $\left( \mathcal{V}_{\beta
_{\alpha },\alpha },s_{\alpha }\right) $, \textit{as viewed from} $\mathcal{M%
}$ (so the available parameters for the definitions come from $M_{\alpha
}\cup \{s_{\alpha }\}$ and consequently $M_{\alpha }\cup \{s_{\alpha
}\}\subseteq N_{\alpha }$)$.$ By Theorem 1.2 we may assume that for some
formula $W(x,y,\overline{m})$ the sentence \textquotedblleft $W$ is a global
well-ordering" is equivalent to a $\Pi _{3}$-statement in\textit{\ }$%
\mathcal{M}$. Therefore, since $n\geq 3$, the statement \textquotedblleft
there is a well-ordering of $\mathrm{V}_{\beta _{\alpha }}$ that is
definable in $(\mathrm{V}_{\beta _{\alpha }},\in )$\textquotedblright\ holds
in $\mathcal{M}$, which immediately shows (by Tarski's elementarity test)
that the statement expressing $\mathcal{N}_{\alpha }\prec \mathcal{V}_{\beta
_{\alpha }}$ holds in $\mathcal{M}.$\footnote{%
This is the only part of the proof that takes advantage of the assumption
that $\mathcal{M}$ carries a definable global well-ordering.} It is
important to have in mind that, as viewed from $\mathcal{M}$, each member of
$\mathcal{N}_{\alpha }$ can be written as the denotation $\delta ^{\mathcal{N%
}_{\alpha }}$ of a definable term $\delta =$ $\delta (\overline{m_{\delta }}%
,c)$ for some $m\in M$ in the language $\mathcal{L}$ described above (where $%
c$ is interpreted by $s_{\alpha })$ so $\delta $ might be of nonstandard
length if $\mathcal{M}$ is not $\omega $-standard (here we are taking
advantage of the definability of a sequence-coding function in $\mathcal{M}%
_{\beta _{\alpha }}$ to reduce the number of parameters of a definable term
that come from $M_{\alpha }$ to one).\medskip

\noindent $(iii)$ Given ordinals $\alpha _{1},\alpha _{2}\in \mathrm{Ord}^{%
\mathcal{M}}$ with $\alpha _{1}<\alpha _{2},$ in $\mathcal{M}$ consider:

\begin{center}
$j_{\alpha _{1},\alpha _{2}}:\mathcal{N}_{\alpha _{1}}\rightarrow \mathcal{N}%
_{\alpha _{2}}$, where $j_{\alpha _{1},\alpha _{2}}(\delta ^{\mathcal{N}%
\alpha _{1}}):=\delta ^{\mathcal{N}_{\alpha _{2}}}.$
\end{center}

\noindent It is not hard to see that $j_{\alpha _{1},\alpha _{2}}$ is an
\textit{elementary} embedding as viewed from $\mathcal{M}$. This follows
from the following key facts:\medskip

\begin{itemize}
\item \noindent $\left( \mathrm{Th}(\mathcal{V}_{\beta _{\alpha _{1}},\alpha
_{1}},s_{\alpha _{1}})\right) ^{\mathcal{M}}=$ $\left( \mathrm{Th}(\mathcal{V%
}_{\beta _{\alpha _{2}},\alpha _{1}},s_{\alpha _{2}})\right) ^{\mathcal{M}}$%
, whenever $\alpha _{1},\alpha _{2}\in \mathrm{Ord}^{\mathcal{M}}$ with $%
\alpha _{1}<\alpha _{2};$ and\medskip

\item $\mathcal{M}\models \mathcal{N}_{\alpha }\prec \mathcal{V}_{\beta
_{\alpha }}$ for each $\alpha \in \mathrm{Ord}^{\mathcal{M}}$.\medskip
\end{itemize}

\noindent $(iv)$ Hence $\left\langle j_{\alpha _{1},\alpha _{2}}:\alpha
_{1}<\alpha _{2}\in \mathrm{Ord}^{\mathcal{M}}\right\rangle $ is a\textit{\
directed system of elementary embeddings}. The desired $\mathcal{N}$ is the
direct limit of this system. Thus, the elements of $\mathcal{N}$ are
equivalence classes $[f]$ of \textquotedblleft strings\textquotedblright\ $f$
of the form:

\begin{center}
\ $f:\{\alpha \in \mathrm{Ord}^{\mathcal{M}}:\alpha \geq \alpha
_{0}\}\rightarrow \bigcup\limits_{\alpha \in \mathrm{Ord}^{\mathcal{M}%
}}N_{\alpha },$
\end{center}

\noindent where $\alpha _{0}\in \mathrm{Ord}^{\mathcal{M}}$ and there is
some $\mathcal{L}$-term $\delta $ such that $m_{\delta }\in M_{\alpha _{0}}$
and $f(\alpha )=\delta ^{\mathcal{N}_{\alpha }}\in N_{\alpha }$ (two strings
are identified iff they agree on a tail of $\mathrm{Ord}^{\mathcal{M}})$. In
particular, for each $\alpha \in \mathrm{Ord}^{\mathcal{M}}$ there is an
embedding:

\begin{center}
$j_{\alpha ,\infty }:\mathcal{N}_{\alpha }\rightarrow \mathcal{N}$, where $%
j_{\alpha ,\infty }(\delta ^{\mathcal{N}_{\alpha }}):=[h],$ and \medskip

$h(\alpha ):=\delta ^{\mathcal{N}_{\alpha }}$ for all $\alpha $ such that $%
m_{\delta }\in M_{\alpha }.$
\end{center}

\noindent A routine variant of Tarski's elementary chains theorem guarantees
that $j_{\alpha ,\infty }$ is an \textit{elementary embedding} for all $%
\alpha \in \mathrm{Ord}^{\mathcal{M}}.\medskip $

\noindent $(v)$ For $m\in M,$ let $f_{m}(\alpha ):=m=\overline{m}^{^{%
\mathcal{N}_{\alpha }}}$ for all $\alpha \in \mathrm{Ord}^{\mathcal{M}}$
such that $m\in M_{\alpha },$ and consider the embedding

\begin{center}
$j:\mathcal{M}\rightarrow \mathcal{N}$, where $j(m):=[f_{m}(\alpha )]$.
\end{center}

\noindent By identifying $m$ with $[f_{m}]$ we can, without loss of
generality, construe $\mathcal{M}$ as a \textit{submodel of} $\mathcal{N}$.$%
\medskip $

A distinguished element of $\mathcal{N}$ is $[g]$, where $g(\alpha
)=s_{\alpha }$ for $\alpha \in \mathrm{Ord}^{\mathcal{M}}.$ $[g]\neq \lbrack
f_{m}]$ for all $m\in M$ since $s_{\alpha }\notin \mathrm{V}_{\alpha }$ for
all $\alpha $ and therefore $g$ and $f_{m}$ differ on a tail of $\alpha \in
\mathrm{Ord}^{\mathcal{M}}.$ This shows that $\mathcal{M}$ is a \textit{%
proper submodel }of $\mathcal{N}$. To see that $\mathcal{N}$ \textit{end}
extends $\mathcal{M}$, suppose $m\in M$ and for some $\mathcal{L}$-definable
term $\delta $, $\delta ^{\mathcal{N}_{\alpha }}\in \overline{m}$ holds in $%
\mathcal{N}_{\alpha }$ for sufficiently large $\alpha $, i.e., for any $%
\alpha $ such that $\{m,m_{\delta }\}\subseteq M_{\alpha }.$ Therefore there
is some $m_{0}\in \mathrm{V}_{\alpha }^{\mathcal{M}}$ such that $\delta ^{%
\mathcal{N}_{\alpha }}=\overline{m_{0}}$ holds in $\mathcal{N}_{\alpha }$
for sufficiently large $\alpha ,$ and therefore also in $\mathcal{N}$, hence
$\mathcal{N}$ end extends $\mathcal{M}$. \medskip

Finally, let's verify that $\mathcal{M}\prec _{\Sigma _{n}}\mathcal{N}$.
Suppose $\mathcal{M}\models \varphi (\overline{m})$, where $\varphi $ is $%
\Sigma _{n}$ and $m\in M.$ Then $\varphi (\overline{m})$ holds for all
sufficiently large $\mathcal{N}_{\alpha }$, since by design we have:

\begin{center}
$\mathcal{N}_{\alpha }\prec \mathcal{M}_{\beta _{\alpha }}\prec _{\Sigma
_{n}}\mathcal{M}$.
\end{center}

\noindent This shows that $\mathcal{N}\models \varphi (\overline{m})$ since,
as observed earlier, each $\mathcal{N}_{\alpha }$ is elementarily embeddable
in $\mathcal{N}$ via $j_{\alpha ,\infty }$.\hfill $\square $\medskip

We are now ready to verify that the tree property for $\mathrm{Ord}$ fails
in the sense of $\mathcal{M}$ for all $\mathcal{M}\models \mathrm{ZFC}$%
.\medskip

\noindent \textbf{2.6.~Theorem.~}\textit{Every model }$\mathcal{M}$ \textit{%
of} $\mathrm{ZFC}$ \textit{carries an} $\mathcal{M}$-\textit{definable} $%
\mathrm{Ord}^{\mathcal{M}}$-\textit{tree no cofinal branch of} \textit{which}
\textit{is} $\mathcal{M}$-\textit{definable}. \medskip

\noindent \textbf{Proof. }The proof splits into two cases, depending on
whether $\mathcal{M}$ satisfies $\exists p\left( \mathrm{V}=\mathrm{HOD}%
(p)\right) $ or not.\footnote{%
Easton proved (in his unpublished dissertation \cite{Easton Thesis}) that
assuming Con($\mathrm{ZF}$) there is a model $\mathcal{M}$ of $\mathrm{ZFC}$
which carries no $\mathcal{M}$-definable global choice function for the
class of pairs in $\mathcal{M}$; and in particular $\exists p\left( \mathrm{V%
}=\mathrm{HOD}(p)\right) $ fails in $\mathcal{M}$. Easton's theorem was
exposited by Felgner \cite[p.231]{Felgner}; for a more recent and
streamlined account, see Hamkins' MathOverflow answer \cite{Joel-failure of
class choice for pairs}.}\medskip

\noindent \textbf{Case 1}. Suppose that $\exists p\left( \mathrm{V}=\mathrm{%
HOD}(p)\right) $ fails in $\mathcal{M}$. Within $\mathrm{ZFC}$ we can define
the tree $\tau _{\mathrm{Choice}}$ whose nodes at level $\alpha $ are choice
functions $f$ for $\mathrm{V}_{\alpha }$, i.e., $f:\mathrm{V}_{\alpha
}\rightarrow \mathrm{V}_{\alpha }$, where $f(x)\in x$ for all nonempty $x\in
\mathrm{V}_{\alpha }$, and the tree ordering is set inclusion. Clearly $%
\mathrm{ZFC}$ can verify that $\tau $ is an $\mathrm{Ord}$-tree. It is also
clear that every $\mathcal{M}$-definable branch of $\tau ^{\mathcal{M}}$ (if
any) is an $\mathcal{M}$-definable global choice function. By Theorem 1.2
this shows that no branch of $\tau _{\mathrm{Choice}}$ is $\mathcal{M}$%
-definable. \medskip

\noindent \textbf{Case 2}. Now suppose $\exists p\left( \mathrm{V}=\mathrm{%
HOD}(p)\right) $ holds in $\mathcal{M}$. Then by Theorem 1.2 there is some $%
\Sigma _{2}$-formula $W(x,y)$ that defines a global well-ordering of $%
\mathcal{M}$. Note that \textquotedblleft $W$ is a global well-ordering" is $%
\Pi _{3}$\textit{-}expressible in $\mathcal{M}$. We claim that for any fixed
$n\geq 3$, no branch of $\tau _{n}^{\mathcal{M}}$ is $\mathcal{M}$%
-definable. If not, then by Lemma 2.4 there is an $\mathcal{M}$-definable
structure $\mathcal{N}$, and an $\mathcal{M}$-definable embedding $j:%
\mathcal{M}\rightarrow \mathcal{N}$ such that $\mathcal{N}$ is a proper is a
$\Sigma _{n}$-e.e.e. of $j(\mathcal{M)}$, which contradicts Theorem
2.1.\hfill $\square $

\bigskip

\begin{center}
\textbf{3.~Consequences of the failure of the definable tree property for
the class of ordinals}\bigskip
\end{center}

\noindent In this section we use Theorem 2.6 to establish further results
about definable combinatorial properties of proper classes within $\mathrm{%
ZFC}$. Our first result improves Theorem 2.6 by combining its proof with
appropriate combinatorial and coding techniques so as to obtain the
description of a single subtree of $^{<\mathrm{Ord}}2$ that is \textrm{Ord}%
-Aronszajn across all models of $\mathrm{ZFC}$; here

\begin{center}
$^{<\mathrm{Ord}}2=\bigcup\limits_{\alpha \in \mathrm{Ord}}{}^{\alpha }2,$
\end{center}

\noindent where${}\ ^{\alpha }2$ is the set of binary sequences of length $%
\alpha .$ The ordering on $2^{<\mathrm{Ord}}$ is `end extension', denoted $%
\sqsubseteq $. Given a tree $\tau $ we say $\tau $ is a \textit{subtree} of $%
\left( ^{<\mathrm{Ord}}2,\ \sqsubseteq \right) $ if each node of $\tau $ is
an element of $^{<\mathrm{Ord}}2$, and the nodes of $\tau $ are ordered by $%
\sqsubseteq $\medskip\

\noindent \textbf{3.1.~Theorem.~}\textit{There is a definable class }$\sigma
$ \textit{that satisfies the following three properties}:\medskip

\noindent \textbf{(a)} $\mathrm{ZFC}\vdash \sigma $ \textit{is a subtree of}
$\left( ^{<\mathrm{Ord}}2,\ \sqsubseteq \right) .$\medskip

\noindent \textbf{(b) }$\mathrm{ZFC}\vdash $ $\sigma $ \textit{is an} $%
\mathrm{Ord}$-\textit{tree}.\medskip

\noindent \textbf{(c)} \textit{For all formulae }$\beta (x,y)$\ \textit{of
set theory,} $\mathrm{ZFC}\vdash $ \textquotedblleft $\{x:\beta (x,y)\}$
\textit{is not a branch of} $\sigma $ \textit{for any parameter }$y$%
\textquotedblright .\medskip

\noindent \textbf{Proof.~}The proof has two stages. In the first stage we
construct an $\mathrm{Ord}$-tree that satisfies properties (b) and (c); and
then in the second stage we construct an appropriate variant of the tree
constructed in the first stage which satisfies properties (a), (b) and (c).
\medskip

\textbf{Stage 1.~}Given $\mathrm{Ord}$-trees $\sigma _{1}=(S_{1},<_{1})$ and
$\sigma _{2}=(S_{2},<_{2})$, let $\sigma _{1}\otimes \sigma _{2}$ be the
tree whose set of nodes is:

\begin{center}
$S_{1}\otimes S_{2}:=\{(p,q)\in S_{1}\times S_{2}:h_{1}(p)=h_{2}(q)\},$
\end{center}

\noindent where $h_{i}(x)$ is the height (level) of $x$, i.e., the ordinal
that measures the order-type of the set of predecessors of $x$ in $\tau _{i}$%
. The ordering on $\sigma _{1}\otimes \sigma _{2}$ is given by:

\begin{center}
$(p,q)\vartriangleleft (p^{\prime },q^{\prime })$ iff $p<_{1}p^{\prime }$
and $q<_{2}q^{\prime }.$
\end{center}

\noindent Routine considerations show that the following two assertions are
verifiable in $\mathrm{ZFC}$:\medskip

\noindent $(i)$ $\sigma _{1}\otimes \sigma _{2}$ is an $\mathrm{Ord}$%
-tree.\medskip

\noindent $(ii)$ Every branch $B$ of $\sigma _{1}\otimes \sigma _{2}$ is of
the form:

\begin{center}
$\{(p,q)\in S_{1}\otimes S_{2}:p\in B_{1}$ and $q\in B_{2}\}$,
\end{center}

\noindent where $B_{i}$ is the branch of $\tau _{i}$ obtained by projecting $%
B$ on its $i$-th coordinate. In particular, for any model $\mathcal{M}%
\models \mathrm{ZFC}$ we have:\medskip

\noindent $(iii)$ If $\left( \sigma _{1}\otimes \sigma _{2}\right) ^{%
\mathcal{M}}$ has an $\mathcal{M}$-definable branch, so do $\sigma _{1}^{%
\mathcal{M}}$ and $\sigma _{2}^{\mathcal{M}}.$\medskip

Let $\sigma _{0}:=\tau _{\mathrm{Choice}}\otimes \tau _{3};$ where $\tau _{%
\mathrm{Choice}}$ and $\tau _{3}$ are as in the proof of Theorem 2.6. It is
easy to see that $\sigma _{0}$ is an \textrm{Ord}-tree (provably in \textrm{%
ZFC}). The proof of Theorem 2.6, coupled with (\textit{iii}) above shows
that no branch of $\sigma _{0}^{\mathcal{M}}$ is $\mathcal{M}$-definable for
any $\mathcal{M}\models \mathrm{ZFC}$. \medskip

\textbf{Stage 2.~}The tools of this stage of the construction are Lemmas
3.1.1 and 3.1.2. Recall that the ordering on both trees $\tau _{\mathrm{%
Choice}}$ and $\tau _{3}$ is set-inclusion $\subseteq .$ \medskip

\noindent \textbf{Lemma 3.1.1.~}\textit{Given} $\mathcal{M}\models \mathrm{%
ZFC}$ \textit{and }$\mathrm{Ord}$\textit{-trees} $\sigma _{1}$ \textit{and} $%
\sigma _{2}$ \textit{in} $\mathcal{M}$ \textit{whose ordering }(\textit{as
viewed in }$\mathcal{M})$ \textit{are set-inclusion, there is an }$\mathcal{M%
}$\textit{-definable }$\mathrm{Ord}$\textit{-tree }$\sigma _{1}\oplus \sigma
_{2}$ \textit{whose ordering is also set-inclusion such that} $\sigma
_{1}\otimes \sigma _{2}$\textit{\ is isomorphic to }$\sigma _{1}\oplus
\sigma _{2}$ \textit{via an }$\mathcal{M}$-\textit{definable isomorphism.}%
\medskip

\noindent \textbf{Proof.} Let $S_{i}$ be the collection of nodes of $\sigma
_{i},$ and consider the tree $\sigma _{1}\oplus \sigma _{2}$ whose sets of
nodes, $S_{1}\oplus S_{2}$, is defined as:

\begin{center}
$\left\{ \left( p\times \{0\}\right) \cup \left( q\times \{1\}\right)
:\left( p,q\right) \in S_{1}\otimes S_{2}\right\} ,$
\end{center}

\noindent and whose ordering is set inclusion. It is easy to see the desired
isomorphism between $\sigma _{1}\otimes \sigma _{2}$ and $\sigma _{1}\oplus
\sigma _{2}$ is described by:

\begin{center}
$(p,q)\mapsto \left( p\times \{0\}\right) \cup \left( q\times \{1\}\right) .$

\hfill $\square $ (Lemma 3.1.1)\medskip
\end{center}

\noindent \textbf{Lemma 3.1.2.~}\textit{Given} $\mathcal{M}\models \mathrm{%
ZFC}$ \textit{and any} $\mathrm{Ord}$\textit{-tree} $\tau $ \textit{in} $%
\mathcal{M}$ \textit{whose ordering }(\textit{as viewed in }$\mathcal{M})$
\textit{is set-inclusion, there is an} $\mathrm{Ord}$\textit{-tree} $%
\widetilde{\tau }$ \textit{of} $\mathcal{M}$ \textit{satisfying the
following properties}:\medskip

\noindent \textbf{(a)} $\widetilde{\tau }=(\widetilde{T},\sqsubseteq ),$%
\textit{\ for some }$\widetilde{T}\subseteq \ ^{<\mathrm{Ord}}2$ .\medskip

\noindent \textbf{(b)} \textit{If} $\widetilde{\tau }$ \textit{has an} $%
\mathcal{M}$-\textit{definable branch}, \textit{then }$\tau $ \textit{has an}
$\mathcal{M}$\textit{-definable branch}.$\medskip $

\noindent \textbf{Proof. }We will first describe a $\mathrm{ZFC}$%
-construction that should be understood to be carried out within $\mathcal{M}
$. Given a set $s$, let $\overline{s}$ be the transitive closure of $\left\{
s\right\} $, and let $\kappa _{s}:=\left\vert \overline{s}\right\vert .$ It
is well-known that given a bijection $g:\overline{s}\rightarrow \left\vert
\overline{s}\right\vert $, $s$ can be canonically coded by some binary
sequence $v_{g}(s)\in \ ^{\left\vert \overline{s}\right\vert }2$. More
specifically, the $\in $ relation on $\overline{s}$ can be readily copied
over $\kappa _{s}$ with the help of $g$ so as to obtain a binary relation $%
R_{g}(s)$ such that $\left( \overline{s},\in \right) \cong \left( \kappa
_{s},R_{g}(s)\right) .$ Since, $R_{g}(s)$ is an extensional well-founded
relation, $s$ can thus be recovered from $R_{g}(s)$ as \textquotedblleft the
top element of the transitive collapse of $R_{g}(s)$\textquotedblright $.$
On the other hand, $R_{g}(s)$ can be coded-up as $X_{g}(s)\subseteq \kappa
_{s}$ with the help of a canonical pairing function $p:\mathrm{Ord}%
^{2}\rightarrow \mathrm{Ord}$. Thus, if $v_{g}(s):$ $\kappa _{s}\rightarrow
\{0,1\}$ is defined as the characteristic function of $X_{g}(s)$, then $%
s=F(v_{g}(s))$, where $F(x)$ is the parameter-free definable class function
given by:

\begin{center}
If $x\in \ ^{<\mathrm{Ord}}2$, and $\overset{x^{\circ }}{\overbrace{%
\{p^{-1}(t):x(t)=1\}}}$ is well-founded, extensional, and has a top element,
then $F(x)$ is the top element of the transitive collapse of $x^{\circ };$
otherwise $F(x)=0.\medskip $
\end{center}

\noindent Given an $\mathrm{Ord}$-tree $\tau =(T,\subseteq )$, let $%
T_{\alpha }$ be the set of elements of $T$ of height $\alpha \in \mathrm{Ord}%
,$ and for $s\in T_{\alpha }$, and $\beta \leq \alpha ,$ let $s_{\beta }$ be
the unique element in $T_{\beta }$ that is a subset of $s.$ Let

\begin{center}
$h_{g}(s):=\bigoplus\limits_{\beta \leq \alpha }v_{g}(s_{\beta }),$
\end{center}

\noindent where $g:\overline{s}\rightarrow \left\vert \overline{s}%
\right\vert $ is a bijection and the operation $\oplus $ is defined as
follows: given a transfinite sequence $\left\langle m_{\beta }:\beta \leq
\alpha \right\rangle $ of binary sequences, $\bigoplus\limits_{\beta \leq
\alpha }m_{\beta }$ is the \textit{ternary} sequence obtained by
concatenating the sequence of sequences $\left\langle m_{\beta }\ast
\left\langle 2\right\rangle :\beta \leq \alpha \right\rangle $, where $%
m_{\beta }\ast \left\langle 2\right\rangle $ is the concatenation of the
sequence $m_{\beta }$ and the sequence $\left\langle 2\right\rangle .$ Thus
the `maximal binary blocks' of $\bigoplus\limits_{\beta \leq \alpha
}m_{\beta }$ are precisely sequences of the form $m_{\beta }$ for some $%
\beta \leq \alpha $. This makes it clear that $s$ can be readily `read off' $%
h_{g}(s)$ as the result of applying $F$ to last binary block of $h_{g}(s)$.

\medskip

Let $\widetilde{T}_{0}:=\{h_{g}(s):s\in T$, and $g$ is a bijection between $%
\overline{s}$ and $\left\vert \overline{s}\right\vert \}.$ We are now ready
to define the desired $\widetilde{T}.$ Fix a canonical embedding $G$ of $^{<%
\mathrm{Ord}}3$ into $^{<\mathrm{Ord}}2$, and let:

\begin{center}
$\widetilde{T}:=\{G(v):v\in \widetilde{T}_{0}\}.$ \medskip
\end{center}

\noindent It is easy to see, using the assumption that $(T,\subseteq )$ is
an $\mathrm{Ord}$-tree, that $\widetilde{\tau }:=(\widetilde{T},\sqsubseteq
) $ is an $\mathrm{Ord}$-tree. Since $\widetilde{T}\subseteq \ ^{<\mathrm{Ord%
}}2$, it remains to show that if $\widetilde{\tau }$ has an $\mathcal{M}$%
-definable branch, then $\tau $ also has an $\mathcal{M}$-definable branch.
Suppose $\widetilde{B}=\{\widetilde{b}_{\alpha }:\alpha \in \mathrm{Ord}^{%
\mathcal{M}}\}$ is a branch of $\widetilde{\tau }.$ Let

\begin{center}
$\widetilde{B}_{0}:=\{G^{-1}(\widetilde{b}_{\alpha }):\alpha \in \mathrm{Ord}%
^{\mathcal{M}}\}$
\end{center}

\noindent Note that $\widetilde{B}_{0}$ is a cofinal branch of the tree $%
\widetilde{\tau }_{0}$; and the maximal binary blocks of $\widetilde{B}_{0}$
form a proper class, and are linearly ordered by set-inclusion (in the sense
of $\mathcal{M}$) by design. Let $B$ be the collection of elements $b\in T$
that are of the form $F(m)$, where $m$ is the last binary block of $G^{-1}(%
\widetilde{b}_{\alpha }).$ Then $B$ is a cofinal branch of $\tau $ and is
definable from $\widetilde{B}.$\hfill $\square $ (Lemma 3.1.2)\medskip

Let $\delta :=\left( \tau _{\mathrm{Choice}}\oplus \tau _{3}\right) $, and $%
\tau :=\widetilde{\delta }$. Theorem 2.6 together with Lemmas 3.1.1 and
3.1.2 make it clear that in every model $\mathcal{M}$ of $\mathrm{ZFC}$, $%
\tau ^{\mathcal{M}}$ is a definably $\mathrm{Ord}$-Aronszajn subtree of $%
\left( ^{<\mathrm{Ord}}2\right) ^{\mathcal{M}}$; so by the completeness
theorem of first order logic, the proof is complete. \hfill $\square $
(Theorem 3.1)$\medskip $

Theorem 3.1 has the following immediate consequence for spartan models of $%
\mathrm{GB+AC}$, where $\mathrm{AC}$ is the axiom of choice for sets:\medskip

\noindent \textbf{3.2.~Corollary.~}\textit{There is a definable class }$%
\sigma $ \textit{in the language of class theory satisfying the following
properties:\medskip }

\noindent \textbf{(a)} $\mathrm{GB+AC}\vdash \sigma $ \textit{is a subtree
of }$^{^{<\mathrm{Ord}}}2$ \textit{and} $\sigma $ \textit{is a proper
class.\medskip }

\noindent \textbf{(b)} \textit{The statement }\textquotedblleft $\sigma $
\textit{is an }$\mathrm{Ord}$-\textit{Aronszajn\ tree}\textquotedblright\
\textit{holds in every spartan model of} $\mathrm{GB+AC}$.\medskip

\noindent \textbf{3.3.~Remark.} It is known \cite[Corollary 2.2.1]{Ali NFUA}
that the set-theoretical consequences of $\mathrm{GB+AC}$ +
\textquotedblleft $\mathrm{Ord}$ has the tree property\textquotedblright\ is
precisely $\mathrm{ZFC}+\Phi $, where $\Phi $ is the scheme whose instances
are of the form \textquotedblleft there is an $n$-Mahlo cardinal $\kappa $
such that $\kappa $ is $n$-correct\textquotedblright , and $n$ ranges over
meta-theoretic natural numbers. Also note that one can derive\ global choice
from local choice in $\mathrm{GB+AC}$ + \textquotedblleft $\mathrm{Ord}$ is
weakly compact\textquotedblright\ (using $\tau _{\mathrm{Choice}}$ of the
proof of Theorem 2.6). Moreover, by an unpublished result of the
first-named-author, there are (non $\omega $-) models $(\mathcal{M},\mathcal{%
S})$ of $\mathrm{GB+AC}$ + \textquotedblleft $\mathrm{Ord}$ has the tree
property\textquotedblright\ in which the partition property $\mathrm{Ord}%
\rightarrow \left( \mathrm{Ord}\right) _{2}^{k}$ fails for some nonstandard $%
k\in \omega ^{\mathcal{M}}$, which implies that for models of $\mathrm{GB+AC}
$, the condition \textquotedblleft $\forall k\in \omega $ $\mathrm{Ord}%
\rightarrow \left( \mathrm{Ord}\right) _{2}^{k}$\textquotedblright\ is
strictly stronger than \textquotedblleft $\mathrm{Ord}$ has the tree
property\textquotedblright .\footnote{%
A similar phenomena occurs in the arithmetic setting in relation to Ramsey's
Theorem: even though the predicative extension $\mathrm{ACA}_{\mathrm{0}}$
of $\mathrm{PA}$ can prove every instance of Ramsey's Theorem of the form $%
\mathbf{\omega }\rightarrow \left( \mathbf{\omega }\right) _{2}^{n}$, where $%
n$ is any meta-theoretic natural number (by a routine arithmetization of any
of the usual proofs of Ramsey's theorem), $\mathrm{ACA}_{\mathrm{0}}$ cannot
prove the stronger statement $\forall k\in \omega $ $\mathbf{\omega }%
\rightarrow \left( \mathbf{\omega }\right) _{2}^{k}$. This natural
incompleteness phenomena follows from a subtle recursion-theoretic theorem
of Jockusch \cite{Jockusch}, which states that for each natural number $%
n\geq 2$ there is a recursive partition $P_{n}$ of $[\omega ]^{n}$ into two
parts such that $P_{n}$ has no infinite $\Sigma _{n}^{0}$-homogeneous
subset. For more detail, see Wang's exposition \cite[p.25]{Wang}; note that
Wang refers to $\mathrm{ACA}_{\mathrm{0}}$ as $\mathrm{PPA}$.}\ But of
course in the Kelley-Morse theory of classes these two statements are
equivalent.\medskip

\noindent \textbf{3.4.~Theorem.} \textit{The definable proper class
partition property fails in every model of }$\mathit{\mathrm{ZFC}}$.\textit{%
\textrm{\ }That is, there is a definable }$2$-\textit{coloring of pairs of
sets having no definable monochromatic proper class. }$\medskip $

\noindent \textbf{Proof}. Let $\tau =(T,\sqsubseteq )$ be as in Theorem 3.1
and $\mathcal{M}\models \mathrm{ZFC}$. We argue in $\mathcal{M}$. For $p,q$
in $T$, we will say that $p$ is \textit{to the right} of $q$, written $%
p\vartriangleright q$, if $p>_{T}q$, or at the point of first difference,
the bit of $p$ is larger than $q$ at that coordinate. Also, as in the proof
of Theorem 3.1, we use $h(p)$ for the height of $p$ in $\tau $. Define a
coloring $f:[T]^{2}\rightarrow \{0,1\}$ by:

\begin{center}
$f(\{p,q\})=\begin{cases}
0,\ \mathrm{if\ }h(p)>h(q),\ \mathrm{and}\ p\vartriangleright q\mathrm{;}\\
1,\ \mathrm{otherwise}.\\
\end{cases}$
\end{center}

\noindent Suppose that $H$ is a definable proper subclass of $T$ that is $f$%
-monochromatic. Next, color pairs from $H$ with color blue if they are of
the same height, and red otherwise. Since the collection of elements of $%
\tau $ of a given height are sets, there cannot be a proper subclass colored
blue, and so we can find a subclass of $H$ with all elements on different
levels. So without loss of generality, all elements on $H$ are on different
levels. If the monochromatic value of pairs from $H$ is $0$, then as one
goes up the tree, the nodes in $H$ are always to the right. Let $B$ consist
of the nodes in $\tau $ that are eventually below the nodes of $H$, that is,
$p\in B$ just in case there is some ordinal $\alpha $ such that all nodes in
$H$ above $\alpha $ are above $p$. It is clear that $B$\ is downward closed.
We claim that $B$ is a branch through $\tau $. $B$ is linearly ordered,
since there can be no first point of nonlinearity: if eventually the nodes
of $H$ are above $p\ast 1$, then they cannot be eventually above $p\ast 0$
(where $\ast $ is the concatenation operation on sequences). Finally, $B$ is
closed under limits, since if $p$ has length $\delta $ and $p|\alpha $ is in
$B$ for all $\alpha <\delta $, then take the supremum of the levels
witnessing that, so you find a single level such that all nodes in $H$ above
that level are above every $p|\alpha $, and so they are above $p$. Thus, $B$
is a branch through $\tau $. But $\tau $ has no definable branches, and so
there cannot be such a monochromatic set $H$. Finally, if the monochromatic
value of $H$ is $1$, then as one goes up, the nodes go to the left, and a
similar argument works. $\square \medskip $

\noindent \textbf{3.5.~Corollary.} $\mathrm{Ord}\rightarrow \left( \mathrm{%
Ord}\right) _{2}^{2}$ \textit{fails for definable classes in} \textit{every
model of }$\mathrm{ZFC}+\exists p\left( \mathrm{V}=\mathrm{HOD}(p)\right) .$
\textit{Indeed, }$\mathrm{Ord}\rightarrow \left( \mathrm{Ord}\right)
_{2}^{2} $ \textit{fails for definable classes in} \textit{every model of }$%
\mathrm{ZFC}$ \textit{in which there is a definable well-ordering of} $^{^{<%
\mathrm{Ord}}}2.\footnote{%
The existence of a global definable well-ordering of $^{^{<\mathrm{Ord}}}2$
is equivalent over $\mathrm{ZF}$ to the so-called Leibniz-Mycielski
principle (LM), explored in \cite{Ali LM}, which includes a result of
Solovay that shows that if $\mathrm{ZF}$ is consistent, then there is a
model of $\mathrm{ZF}+\mathrm{LM}$\ in which \textrm{AC} fails (such a
model, a fortiori, does not carry a parametrically definable global
well-ordering). The conjecture that there is a model of $\mathrm{ZFC}+%
\mathrm{LM}$\ which does not carry a parametrically definable global
well-ordering remains open.}\medskip $

\noindent \textbf{3.6.~Remark. }We do not know whether $\mathrm{Ord}%
\rightarrow \left( \mathrm{Ord}\right) _{2}^{2}$ fails for definable classes
in every model of $\mathrm{ZFC.}$ Some of the usual proofs of the infinite
Ramsey theorem use K\"{o}nig's lemma, which is exactly what is going wrong
with our definably $\mathrm{Ord}$-Aronszajn tree; this suggest that perhaps
there is a definable coloring of pairs of ordinals for which there is no
definable monochromatic proper class of ordinals.$\medskip $

\noindent \textbf{3.7.~Theorem.} \textit{The definable compactness property
fails for }$\mathcal{L}_{\mathrm{\infty },\mathrm{\omega }}$ \textit{in
every model }$\mathcal{M}$\textit{\ of }$\mathrm{ZFC}$.$\medskip $

\noindent \textbf{Proof.} Fix a definable $\mathrm{Ord}$-Aronszajn tree $%
\tau =(T,<_{T})$ of $\mathcal{M}$, and let $\mathcal{L}$ be the language
having a constant $\overline{p}$ for every element $p\in T$ and a binary
relation $<$ for the order of $\tau $, together with a new constant $c$. Let
$\Gamma $ be the theory in $\mathcal{M}$ consisting of the atomic diagram of
$\tau $, together with the assertion that $<$ is a tree order and the
assertions of the form:

\begin{center}
$\varphi _{\alpha }:=\bigvee\limits_{p\in T_{\alpha }}(\overline{p}<c)$.
\end{center}

\noindent That is, $\varphi _{\alpha }$ asserts that the new constant $b$
lies above one of the elements on the $\alpha $-th level $T_{\alpha \text{ }%
} $of $\tau $. In $\mathrm{ZFC}$, having `size $\mathrm{Ord}$' is a stronger
property than `proper class', if global choice fails. Nevertheless, we can
organize $\Gamma $ into an equivalent theory of size $\mathrm{Ord}$ as
follows. Instead of taking the whole atomic diagram as separate statements,
which may not be well-orderable, since we can't seem to well-order the nodes
of $\tau $, we instead for each ordinal $\alpha $ let $\sigma _{\alpha }$ be
the conjunction of the set of atomic assertions that hold in the tree up to
level $\alpha $. Recall that the logic $\mathcal{L}_{\mathrm{Ord,\ \omega }}$
allows the formation of conjunctions of any set of assertions, without
needing to put them into any order. Hence $\Gamma $ is defined in $\mathcal{M%
}$ as $\left\{ \sigma _{\alpha }\wedge \varphi _{\alpha }:\alpha \in \mathrm{%
Ord}\right\} $ plus the sentence that expresses that $<$ is a tree order$.$%
\medskip

Every set-sized subtheory of $\Gamma $ mentions only bounded many sentences
of the form $\sigma _{\alpha }\wedge \varphi _{\alpha },$ so we can find a
model in $\mathcal{M}$ of the subtheory by interpreting $c$ as any element
of the tree $\tau $ on a sufficiently high level. But if there is an $%
\mathcal{M}$-definable model of $\Gamma $, then from that model we can
extract the predecessors of the interpretation of the element $c$, and this
will give an $\mathcal{M}$-definable branch through $\tau $, contradicting
that $\tau $ is definably $\mathrm{Ord}$-Aronszajn in $\mathcal{M}$. \hfill $%
\square \medskip $

We close this section with a conjecture. In what follows $\mathcal{D}_{%
\mathcal{M}}$ is the collection of $\mathcal{M}$-definable subsets of $M$,
and \textquotedblleft $\tau $ is a definably $\mathrm{Ord}$-Suslin tree in $%
\mathcal{M}$\textquotedblright\ means that $(\mathcal{M},\mathcal{D}_{%
\mathcal{M}})$ satisfies \textquotedblleft $\tau $ is an\textit{\ }$\mathrm{%
Ord}$-Aronszajn tree and every anti-chain of $\tau $ has cardinality less
than $\mathrm{Ord}$\textquotedblright $.$\medskip

\noindent \textbf{3.8.~Conjecture.~}\textit{Suppose}\textbf{\ }$\mathcal{M}$
\textit{is a model of }$\mathrm{ZFC}+\mathrm{V}=\mathrm{L}$. \textit{Then
there is some} $\tau _{S}\in $ $\mathcal{D}_{\mathcal{M}}$ \textit{such that
}$\tau _{S}$ \textit{is a definably }$\mathrm{Ord}$\textit{-Suslin tree in }$%
\mathcal{M}$.\textit{\ }\medskip

\noindent Let us motivate the above conjecture. By a theorem of Jensen \cite[%
Theorem VII.1.3]{Devlin} , if $\mathrm{V}=\mathrm{L}$ holds, then every
cardinal $\kappa $ that is not weakly compact carries a $\kappa $-Suslin
tree. The relevant case for us of Jensen's proof is when $\kappa $ is a
strongly inaccessible cardinal. Jensen's proof takes advantage of (1) the
existence of a $\kappa $-Aronszajn tree, and (2) the combinatorial principle
\textquotedblleft for some stationary subset set $E$ of $\kappa $, $\square
_{\kappa }(E)$ holds\textquotedblright . We know, by Theorem 2.6, that the
definable version of (1) can be arranged for $\mathrm{Ord}$. On the other
hand, by adapting Jensen's proof to the definable context, the analogue of
(2) might also be true (using the $\mathrm{V}=\mathrm{L}$ assumption) in $(%
\mathcal{M},\mathcal{D}_{\mathcal{M}}).$ The result in the next section
suggests that perhaps the definable version of (2) holds with the assumption
$\mathrm{V}=\mathrm{L}$ weakened to $\exists p\left( \mathrm{V}=\mathrm{HOD}%
(p)\right) .$ This motivates a stronger form of Conjecture 3.8 in which the
assumption that $\mathrm{V}=\mathrm{L}$ holds in $\mathcal{M}$ is weakened
to the $\mathcal{M}$-definability of a global well-ordering of the
universe.\bigskip

\begin{center}
\textbf{4.~The definable version of }$\Diamond _{\mathrm{Ord}}$ \textbf{and
global definable well-orderings}\bigskip
\end{center}

In this section we show that the definable version of $\Diamond _{\mathrm{Ord%
}}$ holds in a model $\mathcal{M}$ of $\mathrm{ZFC}$ iff $\mathcal{M}$
carries a definable well-ordering of the universe. In light of Theorem 1.2
it follows as a consequence that the definable $\Diamond _{\mathrm{Ord}}$,
although seeming to be fundamentally scheme-theoretic, is actually
expressible in the first-order language of set theory as $\exists p\left(
\mathrm{V}=\mathrm{HOD}(p)\right) $.\medskip

In set theory, the diamond principle asserts the existence of a sequence of
objects, of growing size, such that any large object at the end is very
often anticipated by these approximations. In the case of diamond on the
ordinals, what we will have is a definable sequence of $A_{\alpha }\subseteq
\alpha $, such that for any definable class of ordinals $A$ and any
definable class club set $C$, there are ordinals $\theta \in C$ with $A\cap
\theta =A_{\theta }$. This kind of principle typically allows one to
undertake long constructions that will diagonalize against all the large
objects, by considering and reacting to their approximations $A_{\alpha }$.
Since every large object $A$ is often correctly approximated that way, this
enables many such constructions to succeed.\medskip

\noindent \textbf{4.1.~Theorem.~}\textit{For any model} $\mathcal{M}$
\textit{of} $\mathrm{ZFC}$,\textit{\ if there is an }$\mathcal{M}$\textit{%
-definable well-ordering of the universe, then the definable }$\Diamond _{%
\mathrm{Ord}}$ \textit{holds in }$\mathcal{M}$\textit{. }\medskip

\noindent \textbf{Proof.} We argue in $\mathcal{M}$ to establish the theorem
as a theorem scheme; namely, we shall provide a specific definition within $%
\mathcal{M}$ for the sequence $\vec{A}=\left\langle A_{\alpha }:\theta <%
\mathrm{Ord}\right\rangle $, using the same parameter $p$ as the definition
of the global well-order and with a definition of closely related syntactic
complexity, and then prove as a scheme, a separate statement for each $%
\mathcal{M}$-definable class $A\subseteq \mathrm{Ord}$ and class club $%
C\subseteq \mathrm{Ord}$, that there is some $\theta \in C$ with $A\cap
\theta =A_{\alpha }.$ The definitions of the classes $A$ and $C$ may involve
parameters and have arbitrary complexity.\medskip

Let $\vartriangleleft $ be the definable well-ordering of the universe,
definable by a specific formula using some parameter $p$. We define the $%
\Diamond _{\mathrm{Ord}}$-sequence $\vec{A}=\left\langle A_{\alpha }:\theta <%
\mathrm{Ord}\right\rangle $ by transfinite recursion. Suppose that $\vec{A}%
\upharpoonright \theta $ has been defined. We shall let $A_{\theta
}=\varnothing $ unless $\theta $ is a $\beth $-fixed point above the rank of
$p$ and there is a set $A\subseteq \theta $ and a closed unbounded set $%
C\subseteq \theta $, with both $A$ and $C$ definable in the structure $%
\left( \mathrm{V}_{\theta },\in \right) $ (allowing parameters), such that $%
A\cap \theta \neq A_{\alpha }$ for every $\alpha \in C$. In this case, we
choose the least such pair $(A,C)$, minimizing first on the maximum of the
logical complexities of the definitions of $A$ and of $C$, and then
minimizing on the total length of the defining formulas of $A$ and $C$, and
then minimizing on the G\"{o}del codes of those formulas, and finally on the
parameters used in the definitions, using the well-order $\lhd \
\upharpoonright \mathrm{V}_{\theta }$. For this minimal pair, let $A_{\theta
}=A$. This completes the definition of the sequence $\vec{A}=\left\langle
A_{\alpha }:\theta <\mathrm{Ord}\right\rangle $.\medskip

Let us remark on a subtle point, since the meta-mathematical issues loom
large here. The definition of $\vec{A}$ is internal to the model $\mathcal{M}
$, and at stage $\theta $ we ask about subsets of $\theta $ definable in $%
\left( \mathrm{V}_{\theta },\in \right) $, using the truth predicate for
this structure. If we were to run this definition inside an $\omega $%
-nonstandard model $\mathcal{M}$, it could happen that the minimal formula
we get is nonstandard, and in this case, the set $A$ would not actually be
definable by a standard formula. Also, even when $A$ is definable by a
standard formula, it might be paired (with some constants), with a club set $%
C$ that is defined only by a nonstandard formula (and this is why we
minimize on the maximum of the complexities of the definitions of $A$ and $C$
together). So one must give care in the main argument keeping straight the
distinction between the meta-theoretic natural numbers and the internal
natural numbers of the object theory $\mathrm{ZFC}$.\medskip

Let us now prove that the sequence $\vec{A}$ is indeed a $\Diamond _{\mathrm{%
Ord}}$-sequence for $\mathcal{M}$-definable classes. The argument follows in
spirit the classical proof of $\Diamond $ in the constructible universe $%
\mathrm{L}$, subject to the metamathematical issues we mentioned. If the
sequence $\vec{A}$ does not witness the veracity of the definable $\Diamond
_{\mathrm{Ord}}$ in $\mathcal{M}$, then there is some $\mathcal{M}$%
-definable class $A\subseteq \mathrm{Ord}$, defined in $\mathcal{M}$ by a
specific formula $\varphi $ and parameter $z$, and definable club $%
C\subseteq \mathrm{Ord}$, defined by some $\psi $ and parameter $y$, with $%
A\cap \alpha \neq A_{\alpha }$ for every $\alpha \in C$. We may assume
without loss of generality that these formulas are chosen so as to be
minimal in the sense of the construction, so that the maximum of the
complexities of $\varphi $ and $\psi $ are as small as possible, and the
lengths of the formulas, and the G\"{o}del codes and finally the parameters $%
z,y$ are $\vartriangleleft $-minimal, respectively, successively. Let $m$ be
a sufficiently large natural number, larger than the complexity of the
definitions of $\vartriangleleft ,$ $A$, $C$, and large enough so that the
minimality condition we just discussed is expressible by a $\Sigma _{m}$
formula. Let $\theta $ be any $\Sigma _{m}$-correct ordinal above the ranks
of the parameters used in the definitions. It follows that the restrictions $%
\vartriangleleft \ \upharpoonright \mathrm{V}_{\theta }$ and also $A\cap
\theta $ and $C\cap \theta $ are definable in $\left( \mathrm{V}_{\theta
},\in \right) $ by the same definitions and parameters as their counterparts
in $\mathrm{V}$, that $C\cap \theta $ is club in $\theta $, and $A\cap
\theta $ and $C\cap \theta $ form a form a minimal pair using those
definitions $A\cap \alpha \neq \alpha $ for any $\alpha \in C\cap \theta .$
Thus, by the definition of $\vec{A}$, it follows that $A_{\theta }=A\cap
\theta .$ Since $C\cap \theta $ is unbounded in $\theta $ and $C$ is closed,
it follows that $\theta \in C$ , and so $A_{\theta }=A\cap \theta $
contradicts our assumption about $A$ and $C$. So there are no such
counterexample classes, and thus $\vec{A}$ is a $\Diamond _{\mathrm{Ord}}$%
-sequence with respect to $\mathcal{M}$-definable classes, as claimed.\hfill
$\square $\medskip

\noindent \textbf{4.2.~Theorem.~}\textit{The following are equivalent for }$%
\mathcal{M}\models \mathrm{ZFC}$.\medskip

\noindent $\mathbf{(a)}$ $\mathcal{M}$ \textit{carries an }$\mathcal{M}$%
\textit{-definable global well-ordering.}\medskip

\noindent $\mathbf{(b)}$ $\exists p\left( \mathrm{V}=\mathrm{HOD}(p)\right) $
\textit{holds in }$\mathcal{M}$.\medskip

\noindent $\mathbf{(c)}$ \textit{The definable} $\Diamond _{\mathrm{Ord}}$
\textit{holds in }$\mathcal{M}$. \medskip

\noindent \textbf{Proof.~}We will first give the argument, and then in
Remark 4.3 discuss some issues about the formalization, which involves some
subtle issues.\medskip

\noindent $\mathbf{(a)}\Rightarrow \mathbf{(b)}.$ Suppose that $%
\vartriangleleft $ is a global well-ordering that is definable in $\mathcal{M%
}$ from a parameter $p$. In particular in $\mathcal{M}$ every set has a $%
\vartriangleleft $-minimal element. Let us refine this order by defining $%
x\vartriangleleft ^{\prime }y$, just in case $\rho (x)<\rho (y)$ or $\rho
(x)=\rho (y)$ and $x\vartriangleleft y$ (where $\rho $ is the usual
ordinal-valued rank function). The new order is also a well-order, which now
respects rank. In particular, the order $\vartriangleleft ^{\prime }$ is
set-like, and so every object $x$ is the $\theta $-th element with respect
to the $\vartriangleleft ^{\prime }$-order, for some ordinal $\theta $.
Thus, every object is definable in $\mathcal{M}$ from $p$ and an ordinal,
and so $\mathrm{V}=\mathrm{HOD}(p)$ holds in $\mathcal{M}$, as
desired.\medskip

\noindent $\mathbf{(b)}\Rightarrow \mathbf{(a)}.$ If $\mathcal{M}$ satisfies
$\exists p\ \mathrm{V}=\mathrm{HOD(}p\mathrm{)}$, then we have the canonical
well-order of $\mathrm{HOD}$ using parameter $p$, similar to how one shows
that the axiom of choice holds in $\mathrm{HOD}$. Namely, define $%
x\vartriangleleft y$\ if and only if $\rho (x)<\rho (y)$, or the ranks are
the same, but $x$ is definable from $p$ and ordinal parameters in some $%
\mathrm{V}_{\theta }$ with a smaller $\theta $ than $y$ is, or the ranks are
the same and the $\theta $ is the same, but $x$ is definable in that $%
\mathrm{V}_{\theta }$ by a formula with a smaller G\"{o}del code, or with
the same formula but smaller ordinal parameters. It is easy to see that this
is an $\mathcal{M}$-definable well-ordering of the universe.\medskip

\noindent $\mathbf{(a)}\Rightarrow \mathbf{(c)}.$ This is the content of the
Theorem 4.1.\medskip

\noindent $\mathbf{(c)}\Rightarrow \mathbf{(a)}.$ If $\vec{A}$ is an $%
\mathcal{M}$-definable $\Diamond _{\mathrm{Ord}}$-sequence for $\mathcal{M}$%
-definable classes, then it is easy to see that if $A$ is a set of ordinals
in the sense of $\mathcal{M}$, then $A$ must arise as $A_{\theta }$ for
unboundedly many $\theta \in \mathrm{Ord}^{\mathcal{M}}$. As recalled in the
proof of Lemma 3.1.2, in $\mathrm{ZFC}$ every set is coded by a set of
ordinals. So let us define that $x\vartriangleleft y$, just in case $x$ is
coded by a set of ordinals that appears earlier on $\vec{A}$ than any set of
ordinals coding $y$. This is clearly a well-ordering, since the map sending $%
x$ to the ordinal $\theta $ for which codes $x$ is an $\mathrm{Ord}$-ranking
of $\vartriangleleft $. So there is an $\mathcal{M}$-definable well-ordering
of the universe.\hfill $\square $\medskip

\noindent \textbf{4.3.~Remark.~}An observant reader will notice some
meta-mathematical issues concerning Theorem 4.2. The issue is that
statements (a) and (b) are known to be expressible by statements in the
first-order language of set theory, as single statements, but for statement
(c) we have previously expressed it only as a scheme of first-order
statements. So how can they be equivalent? The answer is that the full
scheme-theoretic content of statement (3) follows already from instances in
which the complexity of the definitions of $A$ and $C$ are bounded.
Basically, once one gets the global well-order, then one can construct a $%
\Diamond _{\mathrm{Ord}}$-sequence that works for all definable classes. In
this sense, we may regard the diamond principle $\Diamond _{\mathrm{Ord}}$
for definable classes as not really a scheme of statements, but rather
equivalent to a single first-order assertion.\medskip

Lastly, let us consider the content of Theorem 4.2 in G\"{o}del-Bernays set
theory or Kelley-Morse set theory. Of course, we know that there can be
models of these theories that do not have $\Diamond _{\mathrm{Ord}}$ in the
full second-order sense. For example, it is relatively consistent with $%
\mathrm{ZFC}$ that an inaccessible cardinal $\kappa $ does not have $%
\Diamond _{\kappa }$, and in this case, the structure $\left( \mathrm{V}%
_{\kappa +1},\mathrm{V}_{\kappa },\in \right) $ will satisfy $\mathrm{GBC}$
and even $\mathrm{KMC}$, but it will not satisfy $\Diamond _{\mathrm{Ord}}$
with respect to all classes, even though it has a well-ordering of the
universe (since there is such a well-ordering in $\mathrm{V}_{\kappa +1}$).
But meanwhile, there will be a $\Diamond _{\mathrm{Ord}}$-sequence that
works with respect to classes that are definable from that well-ordering and
parameters, simply by following the construction given in Theorem
4.2.\medskip

\noindent \textbf{4.4. }A minor adaptation of the proof of Theorem 4.1 shows
that if $\mathcal{M}$\ is a model of $\mathrm{ZFC}$ that carries an $%
\mathcal{M}$-definable global well-ordering, then the definable version of $%
\Diamond _{\mathrm{Ord}}(E)$ holds in\textit{\ }$\mathcal{M}$ for any
definably $\mathcal{M}$-stationary $E\subseteq \mathrm{Ord}^{\mathcal{M}}$:
use the same argument, but only define $A_{\alpha }$ for $\alpha \in E;$ and
in the reflection step of the argument use $\theta \in E\cap C.$ Theorem 4.2
can be also accordingly strengthened.\bigskip

\begin{center}
\textbf{5.~The theory of spartan models of GB}\bigskip
\end{center}

Recall from Section 1 that $\mathrm{GB}_{\mathrm{spa}}$ is the collection of
all sentences that hold in all spartan models of $\mathrm{GB}$. As mentioned
earlier, each theorem scheme of Sections 2 through 4 can be readily
reformulated as demonstrating that a certain sentence belongs to $\mathrm{GB}%
_{\mathrm{spa}}$. Note that the purely set-theoretical consequences of $%
\mathrm{GB}_{\mathrm{spa}}$ coincides with the deductive closure of $\mathrm{%
ZF}$; this is an immediate consequence of coupling the completeness theorem
for first order logic with the fact that $(\mathcal{M},\mathcal{D}_{\mathcal{%
M}})$ is a model of $\mathrm{GB}$ whenever $\mathcal{M}$ is a model of $%
\mathrm{ZF}$. A natural question is whether $\mathrm{GB}_{\mathrm{spa}}$ is
computably axiomatizable. The following result provides a strong negative
answer to this question. \medskip

\noindent \textbf{5.1.~Theorem.~}$\mathrm{GB}_{\mathrm{spa}}$ \textit{is }$%
\Pi _{1}^{1}$-\textit{complete.\medskip }

\noindent \textbf{Proof.~}We need to use both the meta-theoretic natural
numbers, which we will denote by $\omega $, and the object-theoretic natural
numbers, which we denote by $\mathbb{N}$. It is not hard to see that $%
\mathrm{GB}_{\mathrm{spa}}$ has\textit{\ }a\textit{\ }$\Pi _{1}^{1}$%
-description. To see this, consider the following predicates, where $%
r,s\subseteq \omega :\medskip $

\noindent (1) $\mathrm{Sat}_{\mathrm{ZF}}(r)$ expresses \textquotedblleft
the structure canonically coded by $r$ is a model of $\mathrm{ZF}$%
\textquotedblright .

\noindent (2) $s=\mathrm{Def}(r)$ expresses \textquotedblleft $\mathrm{Sat}_{%
\mathrm{ZF}}(r)$\ and $s$ codes the collection of $r$-definable subsets of
the domain of discourse of the structure (coded by) $r$\textquotedblright .

\noindent (3) $\mathrm{Sat}((s,r),\varphi )$ expresses \textquotedblleft $s=%
\mathrm{Def}(r),$ $\varphi $ is a sentence of $\mathcal{L}_{\mathrm{GB}}$,
and the $\mathrm{GB}$-model coded by $(r,s)$ satisfies $\varphi $%
\textquotedblright $.\medskip $

\noindent Usual arguments show that each of the above three predicates is $%
\Delta _{1}^{1}$ in the Baire space. In light of the fact that $\Delta
_{1}^{1}$-predicates are closed under Boolean operations, this makes it
clear that $\mathrm{GB}_{\mathrm{spa}}$ is $\Pi _{1}^{1}$, since by the L%
\"{o}wenheim-Skolem theorem, we have:

\begin{center}
$\varphi \in \mathrm{GB}_{\mathrm{spa}}$ iff $\forall r\subseteq \omega \
\forall s\subseteq \omega \ \left( \left( \mathrm{Sat}_{\mathrm{ZF}%
}(r)\wedge s=\mathrm{Def}(r)\right) \rightarrow \mathrm{Sat}((s,r),\varphi
)\right) $
\end{center}

We next show that $\mathrm{GB}_{\mathrm{spa}}$ is $\Pi _{1}^{1}$-complete.
The revelatory idea here is that within $\mathrm{GB}$ one can define -- via
an existential quantification over classes -- a nonempty `cut'\ $\mathrm{I}$
of ambient natural numbers $\mathbb{N}$ (i.e., a nonempty initial segment $%
\mathrm{I}$ of $\mathbb{N}$ that contains 0 and is closed under successors)
such that:\medskip

\noindent $(\ast )$ If $(\mathcal{M},D_{\mathcal{M}})$ is a spartan model of
\textrm{GB}, then $\mathrm{I}^{(\mathcal{M},D_{\mathcal{M}})}\cong \omega $;
i.e., $\mathrm{I}^{(\mathcal{M},D_{\mathcal{M}})}$ has no nonstandard
elements.\medskip

\noindent The cut $\mathrm{I}$ has a simple definition within $\mathrm{GB}$.
In the definition below $F_{n}$ is the collection of set theoretical
formulae of complexity at most $n,$ where `complexity' can be taken as the
number of occurrences of logical symbols (i.e. the Boolean connectives and
the quantifiers)\footnote{%
The idea of defining the cut $\mathrm{I}$ goes back to\textbf{\ }Mostowski
\cite{Mostowski}, who used it to show that the scheme of induction over $%
\mathbb{N}$ is not provable in \textrm{GB}.\medskip}

\begin{center}
$\mathrm{I}:=\{n\in \mathbb{N}:$ there is a proper class $C$ such that $C$
is the satisfaction-predicate for $F_{n}\}$,
\end{center}

\noindent \textit{The relevant insight is that in spartan models of }$%
\mathrm{GB}$\textit{, the only members of }$\mathrm{I}$\textit{\ are the\
standard natural numbers }$\mathbb{\omega }$\textit{, thanks to Tarski's
undefinability of truth theorem, }which explains the veracity of $(\ast )$.
\medskip

Using $(\ast )$, and the fact that every real can be included in the
standard system of a model of $\mathrm{ZF}$, we will show that every $\Pi
_{1}^{1}$-subset of $\omega $ is many-one reducible to $\mathrm{GB}_{\mathrm{%
spa}}.$ Suppose $P$ is a $\Pi _{1}^{1}$-subset of $\omega $, and let $\omega
^{\omega }$ be the Baire space. Then by Kleene normal form for $\Pi _{1}^{1}$%
-sets \cite{Hartley}, there is some recursive predicate $R(x,y)$ such that:

\begin{center}
$\forall n\left( n\in P\leftrightarrow \forall F\in \omega ^{\omega }\
\exists m\in \omega \ R(F\upharpoonright m,n)\right) ,$
\end{center}

\noindent where $F\upharpoonright m$ is the canonical code for the finite
set of ordered pairs of the form $\left\langle i,F(i)\right\rangle $ with $%
i<m$. Let $\mathsf{R}$ be the formula that numeralwise represents $R$ in $%
\mathrm{GB}$, and given $n\in \omega $, consider the sentence $\varphi _{n}$
in the language of $\mathrm{GB}$ that expresses:

\begin{center}
$\forall s\left( s\in \mathbb{N}\backslash \mathrm{I}\rightarrow \exists
m\in I\ \mathsf{R}^{I}(\mathrm{F}_{s}\upharpoonright m,n)\right) $,
\end{center}

\noindent where $\mathsf{R}^{\mathrm{I}}$ is the result of restricting all
of the quantifiers of the $\mathsf{R}$ to $\mathrm{I}$, and $\mathrm{F}_{s}$
is the function defined in \textrm{GB} with domain $\mathbb{N}$ such that:

\begin{center}
$\mathrm{GB}\vdash \ $\textquotedblleft $\mathrm{F}_{s}(x)$ is the $x$-th
digit of the binary expansion of $s$\textquotedblright .
\end{center}

\noindent It is evident that $n\mapsto \ulcorner \varphi _{n}\urcorner $ is
a computable function. We claim: \medskip

\noindent $(\ast \ast )$ $\forall n(n\in P\leftrightarrow \varphi _{n}\in
\mathrm{GB}_{\mathrm{spa}}).$\medskip

\noindent The left-to-right direction of $(\ast \ast )$ should be clear. The
right-to-left direction is also easy to see, using the fact (proved by a
simple compactness argument) that for every $F\in \omega ^{\omega }$ there
is a non $\omega $-standard model $\mathcal{M}\models \mathrm{ZF}$ and some
nonstandard $s\in \mathbb{N}^{\mathcal{M}}$ such that the `standard part' of
the $\mathcal{M}$-finite function coded by $s$ agrees with $F$, i.e., $%
\forall m\in \omega ~\mathcal{M}\models \left( \mathrm{F}_{s}\upharpoonright
m=F\upharpoonright m\right) .$\hfill $\square $\medskip

\noindent \textbf{5.2.~Remark.~}The above proof strategy can be used to show
that the following theories are also $\Pi _{1}^{1}$-complete:\medskip

\noindent \textbf{(a)} The theory $\left( \mathrm{ACA}_{0}\right) _{\mathrm{%
spa}}$ of all spartan\footnote{%
Spartan models of $\mathrm{ACA}_{0}$ are of the form $\left( \mathcal{M},D_{%
\mathcal{M}}\right) $, where $\mathcal{M}\models \mathrm{PA}.$} models of $%
\mathrm{ACA}_{0}.$\medskip

\noindent \textbf{(b)} The theory of all models of the form $(\mathcal{M}%
,\omega )$, where $\mathcal{M}$ is a model of $\mathrm{ZF}$ or $\mathrm{PA,}$
and $(\mathcal{M},\omega )$ is the expansion of $\mathcal{M}$ by a new
predicate $\omega $ consisting of all \textit{standard} natural numbers in $%
\mathcal{M}$.\medskip

\noindent \textbf{(c)} The theory of all models the form $(\mathcal{M},%
\mathrm{Sat}_{\mathcal{M}})$, where $\mathcal{M}$ is a model of $\mathrm{ZF}$
or $\mathrm{PA}$, and $\mathrm{Sat}_{\mathcal{M}}$ is the satisfaction
predicate for $\mathcal{M}$.

\begin{center}
\begin{tabular}{ll}
{\small Ali Enayat} & {\small Joel David Hamkins} \\
{\small Department of Philosophy, Linguistics, \& Theory of Science} &
{\small The Graduate Center} \\
{\small University of Gothenburg} & {\small The City University of New York}
\\
{\small Box 200, SE Sweden} & {\small 365 Fifth Ave} \\
{\small E-mail: ali.enayat@gu.se} & {\small New York, NY 10016, USA} \\
{\tiny URL:
http://flov.gu.se/english/about/staff?languageId=100001\&userId=xenaal} &
{\small \hspace{0.4in}and} \\
& {\small College of Staten Island} \\
& {\small The City University of New York} \\
& {\small 2800 Victory Boulevard} \\
& {\small Staten Island, NY 10314, USA} \\
& {\small E-mail: jhamkins@gc.cuny.edu} \\
& {\tiny URL: http://jdh.hamkins.org}%
\end{tabular}
\end{center}

\end{document}